\documentclass[11pt]{article}
\usepackage[matrix,arrow,curve,cmtip]{xy}
\usepackage{amssymb}
\usepackage{latexsym}
\usepackage{theorem}
\usepackage[a4paper,width=140mm,height=220mm]{geometry}

\newenvironment{compactitemize}{
\begin{itemize}
\itemsep=-4pt   
\itemindent=0pt 
\labelsep=5pt   
\vspace{-7pt} } 
{\end{itemize}\vspace{-7pt} }

\def\to{\mbox{$\xymatrix@1@C=5mm{\ar@{->}[r]&}$}}
\def\tto{\mbox{$\xymatrix@1@C=5mm{\ar@{=>}[r]&}$}}
\def\halfcirc{\begin{picture}(0,0)\put(0,3){\oval(2,2)[l]}\end{picture}}

\def\incl{\mbox{$\xymatrix@1@C=5mm{\ar@{->}[r]|<{\halfcirc}&}$}}
\def\bkar{\mbox{$\xymatrix@1@C=5mm{\ar@{->}[l]&}$}}
\def\distsign{\begin{picture}(0,0)\put(0,0){\circle{4}}\end{picture}}
\def\dist{\mbox{$\xymatrix@1@C=5mm{\ar@{->}[r]|{\distsign}&}$}}
\def\bkdist{\mbox{$\xymatrix@1@C=5mm{\ar@{->}[l]|{\distsign}&}$}}
\def\biar{\mbox{$\xymatrix@1@C=5mm{\ar@<1.5mm>[r]\ar@<-0.5mm>[r]&}$}}
\def\bidist{\mbox{$\xymatrix@1@C=5mm{\ar@<1.5mm>[r]|{\distsign}\ar@<-0.5mm>[r]|{\distsign}&}$}}
\def\adjar{\mbox{$\xymatrix@1@C=5mm{\ar@<1.5mm>@{<-}[r]\ar@<-0.5mm>[r]&}$}}
\def\adjdist{\mbox{$\xymatrix@1@C=5mm{\ar@<1.5mm>@{<-}[r]|{\distsign}\ar@<-0.5mm>[r]|{\distsign}&}$}}
\def\iso{\mbox{$\xymatrix@1@C=6mm{\ar@{->}[r]^{\sim}&}$}}
\def\doubiso{\mbox{$\xymatrix@1@C=6mm{\ar@{<->}[r]^{\sim}&}$}}
\def\doubar{\mbox{$\xymatrix@1@C=6mm{\ar@{<->}[r]&}$}}

\font\atip = xycmat12 at 10pt
\font\btip = xycmbt12 at 10pt
\def\tip#1{
\mbox{
\begin{picture}(0,0)
\put(-725,-25){\atip\char #1}
\put(-725,-25){\btip\char #1}
\end{picture} }}
\newsavebox{\westsouthwesthead}
\savebox{\westsouthwesthead}{%
\tip{119}}
\newcommand{\wswhead}{\usebox{\westsouthwesthead}}
\def\Endoar#1{
\setlength{\unitlength}{0.01pt}
\ifinner
\mbox{
\begin{picture}(300,1200)(1300,0)
\put(1450,680){\mbox{\footnotesize{$#1$}}}
\put(900,770){\oval(900,900)[t]}
\put(900,770){\oval(900,900)[br]}
\put(900,300){\wswhead}
\end{picture}}
\else
\mbox{
\begin{picture}(200,1400)(1600,400)
\put(2100,1300){\mbox{$#1$}}
\put(1300,1300){\oval(1300,1300)[t]}
\put(1300,1300){\oval(1300,1300)[br]}
\put(1300,600){\wswhead}
\end{picture}}
\fi
\setlength{\unitlength}{1pt}}
\def\endoar{
\setlength{\unitlength}{0.01pt}
\ifinner
\mbox{
\begin{picture}(300,1200)(1300,0)
\put(900,770){\oval(900,900)[t]}
\put(900,770){\oval(900,900)[br]}
\put(900,300){\wswhead}
\end{picture}}
\else
\mbox{
\begin{picture}(200,1400)(1600,400)
\put(1300,1300){\oval(1300,1300)[t]}
\put(1300,1300){\oval(1300,1300)[br]}
\put(1300,600){\wswhead}
\end{picture}}
\fi
\setlength{\unitlength}{1pt}}
\def\Endodist#1{
\setlength{\unitlength}{0.01pt}
\ifinner
\mbox{
\begin{picture}(300,1200)(1300,0)
\put(1600,770){\mbox{\footnotesize{$#1$}}}
\put(900,770){\oval(900,900)[t]}
\put(900,770){\oval(900,900)[br]}
\put(1300,970){\circle{400}}
\put(900,300){\wswhead}
\end{picture}}
\else
\mbox{
\begin{picture}(200,1400)(1600,400)
\put(2100,1300){\mbox{$#1$}}
\put(1300,1300){\oval(1300,1300)[t]}
\put(1300,1300){\oval(1300,1300)[br]}
\put(1850,1600){\circle{400}}
\put(1300,600){\wswhead}
\end{picture}}
\fi
\setlength{\unitlength}{1pt}}

\font\atip = xycmat12 at 10pt
\font\btip = xycmbt12 at 10pt

\def\tip#1{
\mbox{
\begin{picture}(0,0)
\put(-725,-25){\atip\char #1}
\put(-725,-25){\btip\char #1}
\end{picture} }}
\newsavebox{\southhead}
\savebox{\southhead}{%
\tip{15}}
\newcommand{\shead}{\usebox{\southhead}}

\newsavebox{\northhead}
\savebox{\northhead}{%
\tip{79}}

\newsavebox{\northnorthnortheasthead}
\savebox{\northnorthnortheasthead}{%
\tip{75}}

\newsavebox{\northnortheasthead}
\savebox{\northnortheasthead}{%
\tip{71}}

\def\down{
\setlength{\unitlength}{0.01pt}
\mbox{
\begin{picture}(300,0)
\put(-100,-200){\line(0,1){1000}}
\put(-100,-200){\shead}
\end{picture}}
\setlength{\unitlength}{1pt}}

\def\ddown{
\setlength{\unitlength}{0.01pt}
\mbox{
\begin{picture}(300,0)
\put(-100,-200){\line(0,1){1000}}
\put(-100,-200){\shead}
\put(-100,0){\shead}
\end{picture}}
\setlength{\unitlength}{1pt}}

\newtheorem{theorem}{Theorem}[section]

\newtheorem{definition}[theorem]{Definition} 
\newtheorem{proposition}[theorem]{Proposition}
\newtheorem{corollary}[theorem]{Corollary}
{\theorembodyfont{\upshape}\newtheorem{example}[theorem]{Example}}
\newcommand{\proof}{\noindent {\it Proof\ }: }
\def\endofproof{$\mbox{ }\hfill\Box$\par\vspace{1.8mm}\noindent}

\def\Mnd{{\sf Mnd}}
\def\Idm{{\sf Idm}}
\def\Y{{\cal Y}}
\def\R{{\cal R}}
\def\RSDist{{\sf RSDist}}
\def\RSCat{{\sf RSCat}}

\def\ol#1{\overline{#1}}
\def\+{^{\dagger}}
\def\etal{~{\it et~al.}}

\def\Cocont{{\sf Cocont}}
\def\Distrib{{\sf Distrib}}

\def\Sh{{\sf Sh}}

\def\:{\colon}
\def\1{{\bf 1}}

\def\2{{\bf 2}}
\def\3{{\bf 3}}

\def\Set{{\sf Set}}
\def\QUANT{{\sf QUANT}}

\def\CAT{{\sf CAT}}

\def\Set{{\sf Set}}

\def\Sup{{\sf Sup}}
\def\Cat{{\sf Cat}}
\def\Matr{{\sf Matr}}

\def\Dist{{\sf Dist}}

\def\Cat{{\sf Cat}}
\def\CAT{{\sf CAT}}
\def\SCat{{\sf SCat}}

\def\Q{{\cal Q}}

\def\C{{\cal C\!\mbox{ }}}

\def\P{{\cal P}}
\def\V{{\cal V}}
\def\W{{\cal W}}

\def\colim{\mathop{\rm colim}}
\def\lim{\mathop{\rm lim}}
\def\bbA{\mathbb{A}}
\def\bbB{\mathbb{B}}
\def\bbC{\mathbb{C}}
\def\bbD{\mathbb{D}}

\def\tensor{\otimes}
\def\<{\langle}
\def\>{\rangle}

\title{Categorical structures enriched in a quantaloid: \\ regular presheaves, regular
semicategories}
\author{Isar Stubbe\footnote{D\'epartement de Math\'ematique, Universit\'e de Louvain, Chemin du Cyclotron 2,
1348 Louvain-la-Neuve (Belgique), {\tt
i.stubbe@math.ucl.ac.be}.}}
\date{June 22, 2004}

\begin{document}

\maketitle

\begin{quote}{\small

{\bf R\'esum\'e.} On \'etudie les pr\'efaisceaux sur des semicat\'egories enrichies dans un quantalo\"{i}de: cela donne lieu \`a la notion de pr\'efaisceau r\'egulier. Une semicat\'egorie est r\'eguli\`ere si tous les pr\'efaisceaux repr\'esentables sont r\'eguliers, et ses pr\'efaisceaux r\'eguliers forment alors une (co)localistation essentielle de la cat\'egorie de tous ses pr\'efaisceaux. La notion de semidistributeur r\'egulier permet d'\'etablir l'\'e\-qui\-va\-len\-ce de Morita des semicat\'egories r\'eguli\`eres. Les ordres continus et les $\Omega$-ensembles fournissent des exemples.\\
{\bf Mots cl\'es:} quantalo\"{i}de, semicat\'egorie, pr\'efaisceau, r\'egularit\'e, \'equivalence de Morita, ordre continu, $\Omega$-ensemble

}\end{quote}

\section{Introduction}

In [Moens\etal, 2002] the theory of regular modules on an $R$-algebra without unit, for $R$ a commutative ring, was generalized to a theory of regular presheaves on a $\V$-enriched semicategory, for $\V$ a symmetric monoidal closed base category. As a monoidal category $\V$ is a one-object bicategory, it is natural to ask in how far in the above the base $\V$ can be replaced by a bicategory $\W$ (thus necessarily loosing symmetry of the tensor). Here we present such a theory of regular presheaves on a $\Q$-enriched semicategory, where now $\Q$ is any (small) quantaloid.
\par
A quantaloid is a $\Sup$-enriched category; it is thus in particular a bicategory. There is a theory of categories enriched in a quantaloid $\Q$, as particular case of categories enriched in a bicategory. A presentation thereof is given in [Stubbe, 2004a] which is our reference for all the basic notions and results concerning $\Q$-categories that we may need further on. A $\Q$-semicategory is then simply a ``$\Q$-category without unit-inequalities''. 
\par
A presheaf on a $\Q$-semicategory $\bbA$ is formally the same thing as a presheaf on the free $\Q$-category on $\bbA$. Thus the presheaves on $\bbA$ constitute a $\Q$-category $\P\bbA$. But the behaviour of those presheaves on a semicategory is radically different from that of presheaves on a category: although there is a Yoneda semifunctor $\bbA\to\P\bbA$, sending an object to a representable presheaf, there is no ``Yoneda lemma'' for semicategories! This leads us to consider two specific kinds of presheaves on a $\Q$-semicategory: those which are canonically the colimit of representable presheaves, that we call regular presheaves, and those for which homming with a representable gives its value in the representing object, that we call Yoneda presheaves. A semicategory for which all representables are regular, is called a regular semicategory (and a semicategory for which all representables are Yoneda, is a category). Our treatment of these matters is sometimes different from [Moens\etal, 2002]; but we recover their results (modulo a translation from $\V$-enrichment to $\Q$-enrichement), so we use their terminology.
\par
An important result is the following: regular presheaves on a regular $\Q$-semi\-cat\-e\-go\-ry $\bbA$ form a $\Q$-category $\R\bbA$ which is an essential (co)localisation of the category $\P\bbA$ of all presheaves on $\bbA$; and the image of the ultimate right adjoint involved in that (co)localisation is the $\Q$-category $\Y\bbA$ of Yoneda presheaves on $\bbA$. From there on, and with a suitable notion of regular semidistributor between regular semicategories, it is a matter of straightforwardly generalizing $\Q$-category theory to obtain an aspect of Morita equivalence for semicategories: for regular $\Q$-semicategories $\bbA$ and $\bbB$, $\R\bbA$ is equivalent to $\R\bbB$ if and only if $\bbA$ and $\bbB$ are isomorphic in the quantaloid of regular semidistributors. That  regular $\Q$-semicategories are profoundly different from $\Q$-categories, is then exemplified by the fact that not every regular $\Q$-semicategory is Morita equivalent to a $\Q$-category.
\par
An inspiration for this work has been the search for a suitable notion of ``sheaf on a quantale (or quantaloid) $\Q$'', generalizing sheaves on a locale $\Omega$. Some have suggested that the regular $\Q$-semicategories should play this r\^ole [Van den Bossche, 1995; van der Plancke, 1997], for an $\Omega$-set is indeed exactly the same thing as a symmetric regular $\Omega$-semicategory\footnote{However those authors did not study ``regularity'' {\em an sich}, as we do here. Instead, based on heuristic arguments they proposed an algebraic gadget called $\Q$-set as generalization of $\Omega$-set. But it turns out that such a $\Q$-set is precisely a regular $\Q$-semicategory.}. This is not entirely our opinion, even though the regular presheaves on an $\Omega$-set are precisely its $\Omega$-subsets. We will present our point of view on ``(ordered) sheaves on a quantaloid'' in [Stubbe, 2004b] for it is too involved to include it here; but the results of the present study -- in particular the notions of regular semicategory, regular semidistributor and regular semifunctor -- are crucial for that further development.
\par
Another interesting point of view on regular semicategories is the following: understanding a $\Q$-semicategory as (a classification into $\Q$ of) a transitive relation on a ($\Q$-typed) set, the regularity of that semicategory comes down to the interpolation property for that transitive relation. For example, a transitive relation $\prec$ on a set $A$ is a $\2$-enriched semicategory; its regularity would mean that for every $a\prec b$ in $A$ there is an $x$ in $A$ such that $a\prec x\prec b$. The way-below relation on a continuous order is thus a regular $\2$-semicategory; moreover the (covariant) regular presheaves on such a regular semicategory are its Scott-open subsets, and the (contravariant) Yoneda presheaves the Scott-closed ones. This gives an understanding of the $\Q$-category $\R\bbA$ of regular presheaves on a regular $\Q$-semicategory as some kind of Scott-topology on the objects of $\bbA$.
\par
Finally we should mention that the quantaloid of regular $\Q$-semicategories and regular distributors enjoys a universal property: it is the Cauchy completion of the base quantaloid $\Q$, i.e.~it is its completion for absolute weighted (co)limits in the (illegitimate) category of quantaloids and quantaloid homomorphisms. As this is not the main theme of this paper, we have merely added some useful hints (with references) in an appendix.

\section{Semicategories, semidistributors and semifunctors} 

Throughout $\Q$ denotes a {\em small} quantaloid. 
\begin{definition}\label{401} A {\em
$\Q$-semicategory} $\bbA$ consists of
\begin{compactitemize}
\item objects: a $\Q$-typed set\footnote{A {\em $\Q$-typed set $X$} is an
object of the slice category $\Set/\Q_0$ of sets over the object-set of $\Q$. In other
terms, such is a set
$X$ to every element of which is associated an object of $\Q$: for every $x\in X$ there
is a $tx$ in $\Q$ (which is called the {\em type} of $x$ in $\Q$). The notation with a
``$t$'' for the types of elements in a $\Q$-typed set is generic: even for two
different $\Q$-typed sets $X$ and $Y$, the type of $x\in X$ is written $tx$ and that of
$y\in Y$ is $ty$.} $\bbA_0$,
\item hom-arrows:  for all $a,a'\in\bbA_0$, a $\Q$-arrow $\bbA(a',a):ta\to ta'$,
\end{compactitemize} satisfying
\begin{compactitemize}
\item composition-inequalities: for all $a,a',a''\in\bbA_0$, 
$\bbA(a'',a')\circ\bbA(a',a)\leq\bbA(a'',a)$ in $\Q$.
\end{compactitemize}
\end{definition} 
\begin{definition}\label{402} For $\Q$-semicategories $\bbA$ and $\bbB$, a {\em
semidistributor} $\Phi\:
\bbA\dist\bbB$ is determined by
\begin{compactitemize}
\item distributor-arrows: for all $a\in\bbA_0, b\in\bbB_0$, a $\Q$-arrow
$\Phi(b,a)\: ta\to tb$ 
\end{compactitemize} satisfying
\begin{compactitemize}
\item action-inequalities: for all $a,a'\in\bbA_0$ and $b,b'\in\bbB_0$, 
$\Phi(b,a')\circ\bbA(a',a)\leq\Phi(b,a)$ and $\bbB(b,b')\circ\Phi(b',a)\leq\Phi(b,a)$.
\end{compactitemize}
\end{definition} 
\begin{definition}\label{403} For $\Q$-semicategories $\bbA$ and
$\bbB$, a {\em semifunctor}
$F\:\bbA\to\bbB$ is 
\begin{compactitemize}
\item object-mapping: a map
$F\:\bbA_0\to\bbB_0\:a\mapsto Fa$
\end{compactitemize} satisfying
\begin{compactitemize}
\item type-equalities: for all $a\in\bbA_0,\ t(Fa)=ta$,
\item action-inequalities: for all $a,a'\in\bbA_0,\ \bbA(a',a)\leq\bbB(Fa',Fa)$.
\end{compactitemize}
\end{definition} 
\par A $\Q$-semicategory is precisely a ``$\Q$-category without
unit-inequalities'': indeed, take the data and axioms for a $\Q$-semicategory $\bbA$ and add the requirement that
\begin{compactitemize}
\item unit-inequalitities: for all $a\in\bbA_0$, $1_{ta}\leq\bbA(a,a)$,
\end{compactitemize}
and one has the definition of ``$\Q$-category''. Consequently, 
any $\Q$-category is a
$\Q$-semicategory (but the converse is obviously false). 
\par
When $\bbA$ and $\bbB$ are $\Q$-categories, then a semidistributor from
$\bbA$ to $\bbB$ is a distributor, and a
semifunctor from $\bbA$ to $\bbB$ is a functor. 
The reason is that neither
a distributor nor a functor between
$\Q$-categories make any reference to the unit-inequalities in its domain and codomain;
this is ultimately a consequence of the base quantaloid $\Q$ being locally ordered. The use of the words ``semidistributor'' and ``semifunctor'' is meant to stress that domain and codomain may be $\Q$-semicategories instead of categories.
\begin{example}\label{ex1}
Let $\2$ denote the two-element Boolean algebra; regard it as a one-object quantaloid. A $\2$-semicategory $\bbA$ is a set $A=\bbA_0$ with a transitive binary relation $\prec$ defined as $a\prec b\iff \bbA(a,b)=1$. So, for example, any strict order is a $\2$-semicategory (but not a $\2$-category). A semifunctor $F\:\bbA\to\bbB$ between $\2$-semicategories is then a relation-preserving map $\bbA_0\to\bbB_0\:x\mapsto Fx$ ($x\prec y$ implies $Fx\prec Fy$); and a semidistributor $\Phi\:\bbA\dist\bbB$ is a relation $\Phi\subseteq\bbB_0\times\bbA_0$ which is ``up-closed'' in $\bbA$ and ``down-closed'' in $\bbB$.
\end{example}
\begin{example}\label{ex2}
The objects of the topos $\Sh(\Omega)$ of sheaves on a locale $\Omega$ may be described as $\Omega$-sets: such is a set $A$ together with a map $[\cdot = \cdot]\:A\times A\to\Omega$ satisfying axioms that say precisely that $\bbA$, with $\bbA_0=A$ and $\bbA(a,b)=[a=b]$, is a symmetric $\Omega$-semicategory, with $\Omega$ viewed as one-object quantaloid. (``Symmetry'' obviously means that $\bbA(a,b)=\bbA(b,a)$ for all $a,b\in\bbA_0$). 
\end{example}
\par  Clearly semifunctors can be composed by composing the object-mappings, and the
identity object-mapping is the identity semifunctor for this composition, precisely as was
the case for functors. 
\begin{proposition}
$\Q$-semicategories and semifunctors form a category
$\SCat(\Q)$, of which $\Cat(\Q)$ is a full subcategory. 
\end{proposition}
Every $\Q$-semicategory $\bbA$ determines a $\Q$-category $\ol{\bbA}$ in
the following way: 
\begin{compactitemize}
\item objects: $(\ol{\bbA})_0=\bbA_0$ as $\Q$-typed sets,
\item hom-arrows: for $a,a'\in(\ol{\bbA})_0$,
$\ol{\bbA}(a',a)=\left\{\begin{array}{ll}
\bbA(a',a)&\mbox{ if }a'\neq a,\\
\bbA(a,a)\bigvee 1_{ta}&\mbox{ if }a'=a.
\end{array}\right.$
\end{compactitemize}
The inclusion semifunctor
$i_{\bbA}\:\bbA\to\ol{\bbA}\:a\mapsto a$ displays $\ol{\bbA}$ as free $\Q$-category on the $\Q$-semicategory $\bbA$, as the following is a trivial fact.
\begin{proposition}\label{401.0}
For a $\Q$-semicategory $\bbA$ and a $\Q$-category $\bbC$, a type-preserving
object-mapping $F\:\bbA_0\to\bbC_0\:A\mapsto Fa$ determines a semifunctor $F\:\bbA\to\bbC$ if and only if it determines a functor $F\:\ol{\bbA}\to\bbC$.
\end{proposition}
Explicitly, for every semicategory $\bbA$ and category $\bbC$, 
\begin{equation}\label{407}
\Cat(\Q)(\ol{\bbA},\bbC)=\SCat(\Q)(\bbA,\bbC), 
\end{equation}  
so there is an adjuction $\ol{(-)}\dashv i\:\Cat(\Q)\adjar\SCat(\Q)$
displaying $\Q$-categories and functors as a full reflective subcategory of
$\Q$-semi\-cat\-e\-go\-ries and semifunctors.
\par  
As for semidistributors, for two $\Q$-semicategories $\bbA$ and $\bbB$, denote
$\Distrib(\bbA,\bbB)$ for the set of semidistributors with domain $\bbA$ and codomain $\bbB$. The proof of the following is straightforward.
\begin{proposition}\label{406} For any two $\Q$-semicategories $\bbA$ and $\bbB$, a
collection of $\Q$-arrows
$\Phi(b,a)\:ta\to tb$, for $a\in\bbA_0$ and $b\in\bbB_0$, determines a
semidistributor $\Phi\:\bbA\dist\bbB$ between semicategories if and only if it determines a
distributor $\Phi\:\ol{\bbA}\dist\ol{\bbB}$ between categories.
\end{proposition}
That is to say, for all $\Q$-semicategories $\bbA$ and $\bbB$,
\begin{equation}\label{408}
\Distrib(\bbA,\bbB)=\Dist(\Q)(\ol{\bbA},\ol{\bbB}).
\end{equation} 
\par
The right hand side of (\ref{408}) is a sup-lattice: namely a
hom-object of the quantaloid $\Dist(\Q)$ of $\Q$-categories and distributors. Therefore the left hand side is a sup-lattice too: for semidistributors $(\Phi_i\:\bbA\dist\bbB)_{i\in I}$ the supremum   $\bigvee_i\Phi_i\:\bbA\dist\bbB$ has elements
\begin{equation}\label{408.1}
(\bigvee_i\Phi_i)(b,a)=\bigvee_i\Phi_i(b,a)\mbox{ for $a\in\bbA_0$ and $b\in\bbB_0$}.
\end{equation}
(That is to say, $\bigvee_i\Phi_i$ is the semidistributor corresponding through (\ref{408}) to the supremum of the distributors $\Phi_i\:\ol{\bbA}\dist\ol{\bbB}$ between the free categories.)
\par
Further, (\ref{408}) also suggests a ``composition law'' for semidistributors:  given semidistributors $\Phi\:\bbA\dist\bbB$ and $\Psi\:\bbB\dist\bbC$, let the ``composite'' $\Psi\tensor\Phi\:\bbA\dist\bbC$ have dis\-tri\-bu\-tor-ar\-rows \begin{equation}\label{408.0}
\Big(\Psi\tensor\Phi\Big)(c,a)=
\bigvee_{b\in\bbB_0}\Psi(c,b)\circ\Phi(b,a)\mbox{ for $a\in\bbA_0$ and $c\in\bbC_0$.}
\end{equation}
(So $\Psi\tensor\Phi$ corresponds through (\ref{408}) to the composition in $\Dist(\Q)$ of $\Phi\:\ol{\bbA}\dist\ol{\bbB}$ with $\Psi\:\ol{\bbB}\dist\ol{\bbC}$.) 
\par
Surely the ``composition'' in (\ref{408.0}) distributes over the suprema in (\ref{408.1}). But the sup-lattices of semidistributors between semicategories are {\em not} the hom-objects of some ``quantaloid of all semicategories and semidistributors'' (and therefore we've put the word ``composition'' in the previous paragraph between inverted commas): even though any $\Q$-semicategory $\bbA$ still determines a
semidistributor $\bbA\:\bbA\dist\bbA$, with distributor-arrows 
\begin{equation}\label{408.2}
\bbA(a',a)\:ta\to ta'\mbox{ for $a,a'\in\bbA_0$},
\end{equation} 
this endo-semidistributor is in general {\em not the identity} for the ``composition'' in (\ref{408.0})! (The definition of `regular semicategory' and `regular semidistributor' fix precisely this problem---see especially section \ref{465}.)

\section{Presheaves on a semicategory}

\subsection*{Contravariant presheaves}

Recall that, given a $\Q$-category $\bbC$, a {\em contravariant presheaf of type $X$ on $\bbC$} is a distributor with domain\footnote{For an object $X$ of $\Q$ we denote $*_X$ for the $\Q$-category with one object and hom-arrow $1_X$. } $*_X$ and codomain $\bbC$. Those presheaves are the objects of a new $\Q$-category $\P\bbC$. The hom-arrow $\P\bbC(\psi,\phi)\:X\to Y$ between objects $\phi\:*_X\dist\bbC$ to  $\psi\:*_Y\dist\bbC$ is (the single element of) the lifting $[\psi,\phi]$ in the quantaloid $\Dist(\Q)$; explicitly, $\P\bbA(\psi,\phi)=\bigwedge_{c\in\bbC}[\psi(c),\phi(c)]$ (an infimum of liftings in $\Q$). All this now motivates the following.
\begin{definition}\label{409} A {\em contravariant presheaf of type $X$} on a $\Q$-semicategory $\bbA$ is a semidistributor from $*_X$ into $\bbA$.
\end{definition}  
\begin{proposition}\label{409.1}
The contravariant presheaves on a $\Q$-semicategory $\bbA$ form a $\Q$-category\footnote{There can be no notational confusion concerning the category of presheaves on a $\Q$-category $\bbC$, because the presheaves on $\bbC$-as-semicategory {\em are} the presheaves on $\bbC$-as-category.} $\P\bbA$ with: 
\begin{compactitemize}
\item objects: the $\Q$-typed set of contravariant presheaves on $\bbA$,
\item hom-arrows: for $\phi\:*_X\dist\bbA$ and $\psi\:*_Y\dist\bbA$, $\P\bbA(\psi,\phi)=\bigwedge_{a\in\bbA_0}[\psi(a),\phi(a)]$.
\end{compactitemize}
This $\Q$-category $\P\bbA$ is equivalent to the $\Q$-category $\P\ol{\bbA}$ of presheaves on the free $\Q$-category $\ol{\bbA}$.
\end{proposition}
\proof
The point is that a presheaf $\phi\:*_X\dist\bbA$ on a $\Q$-semicategory $\bbA$ may by (\ref{408}) be viewed as a presheaf $\phi\:*_X\dist\ol{\bbA}$ on the free category $\ol{\bbA}$, and {\em vice versa}; so the objects of $\P\bbA$ and $\P\ol{\bbA}$ are the same. Viewing two semidistributors $\phi\:*_X\dist\bbA$ and $\psi\:*_Y\dist\bbA$ as distributors $\phi\:*_X\dist\ol{\bbA}$ and $\psi\:*_Y\dist\ol{\bbA}$, the $\Q$-arrows $\P\bbA(\psi,\phi)$ and $\P\ol{\bbA}(\psi,\phi)$ are equal. So $\P\bbA=\P\ol{\bbA}$, which proves in particular that $\P\bbA$ is a $\Q$-category. 
\endofproof
\par
It may thus seem that the notion of ``presheaf on a semicategory'' reduces entirely to that of ``presheaf on a category''. But the notion of presheaf presented in \ref{409} is not the interesting one in the theory of $\Q$-semicategories: the presheaf category $\P\bbA$ says a lot about the free category $\ol{\bbA}$, but very little about the semicategory $\bbA$ (see \ref{ex5.3} for a clear example)! Our reason for introducing the $\Q$-category $\P\bbA$ of {\em all} presheaves on a semicategory $\bbA$, is that (it exists anyway and) it allows us to analyze precisely how a particular essential (co)localization of $\P\bbA$ contains those presheaves which are indeed the pertinent ones in the study of $\Q$-semicategories (see \ref{412}, \ref{418} and \ref{422} further on).
\par
Recall that the $\Q$-category $\P\bbC$ of contravariant presheaves on a $\Q$-category $\bbC$ classifies distributors with codomain $\bbC$: for every $\Q$-category $\bbD$, $\Dist(\Q)(\bbD,\bbC)\cong\Cat(\Q)(\bbD,\P\bbC)$. Explicitly, a distributor $\Phi\:\bbD\dist\bbC$ determines the functor $$F_{\Phi}\:\bbD\to\P\bbC\:d\mapsto\Big(\Phi(-,d)\:*_{td}\dist\bbC\Big).$$ Using
(\ref{407}) and (\ref{408}), and the fact that $\P\bbA=\P\ol{\bbA}$, it follows that for $\Q$-semicategories $\bbA$ and $\bbB$,
$$\Distrib(\bbB,\bbA)=\Dist(\Q)(\ol{\bbB},\ol{\bbA})
\cong\Cat(\Q)(\ol{\bbB},\P\ol{\bbA})=\SCat(\Q)(\bbB,\P\bbA).$$
So the presheaf category $\P\bbA$ classifies
semidistributors with codomain $\bbA$: such a semidistributor $\Phi\:\bbB\dist\bbA$ is sent to the semifunctor $$F_{\Phi}\:\bbB\to\P\bbA\:b\mapsto\Big(\Phi(-,b)\:*_{tb}\dist\bbA\Big).$$
In particular does, for any $\Q$-semicategory $\bbA$, the semidistributor $\bbA\:\bbA\dist\bbA$ induce the semifunctor 
$$Y_{\bbA}\:\bbA\to\P\bbA\:a\mapsto\Big(\bbA(-,a)\:*_{ta}\dist\bbA\Big);$$
it is the {\em Yoneda semifunctor} for $\bbA$. An object in the image of $Y_{\bbA}$ is a {\em representable contravariant presheaf}. 
\par
A word of caution is in order here. By the universal property of the
free category
$\ol{\bbA}$, to the Yoneda semifunctor
$Y_{\bbA}\:\bbA\to\P{\bbA}$ corresponds a functor
$\ol{Y_{\bbA}}\:\ol{\bbA}\to\P{\bbA}$. On the other hand, there is the Yoneda embedding 
$Y_{\ol{\bbA}}\:\ol{\bbA}\to\P{\bbA}$ for the $\Q$-category
$\ol{\bbA}$. The functors $\ol{Y_{\bbA}}$ and $Y_{\ol{\bbA}}$ are very different!
\begin{example}\label{ex3}
Considering a strict order $(A,<)$ as a $\2$-semicategory $\bbA$, the free $\2$-category $\ol{\bbA}$ is the (reflexive) order $(A,\leq)$ defined as $a\leq b$ if and only if $a<b$ or $a=b$. A presheaf on $\bbA$ is by definition a presheaf on $\ol{\bbA}$, that is, a downset of $(A,\leq)$. The Yoneda semifunctor $Y_{\bbA}\:\bbA\to\P\bbA$ sends an element $a\in A$ to the {\em strict} principal downset $\ddown a=\{x\in A\mid x<a\}$, hence so does the functor $\ol{Y_{\bbA}}\:\ol{\bbA}\to\P\bbA$. The Yoneda embedding $Y_{\ol{\bbA}}\:\ol{\bbA}\to\P\bbA$ on the other hand sends an element $a\in A$ to the principal downset $\down a=\{x\in A\mid x\leq a\}$.
\end{example}

\subsection*{Yoneda presheaves}

For presheaves on a
$\Q$-category $\bbC$ the ``Yoneda lemma'' holds: for any $\phi\in(\P\bbC)_0$ and $c\in\bbC_0$,
$Y_{\bbC}\:\bbC\to\P\bbC$ satisfies
$\P\bbC(Y_{\bbC}c,\phi)=\phi(c)$. As a consequence, the Yoneda embedding
$Y_{\bbC}\:\bbC\to\P\bbC$ is a fully faithful functor: $\bbC(c',c)=\P\bbC(Y_{\bbC}c',Y_{\bbC}c)$ for all $c,c'\in\bbC_0$.
\par For a $\Q$-semicategory $\bbA$, the Yoneda semifunctor $Y_{\bbA}\:\bbA\to\P{\bbA}$ cannot be fully faithful unless $\bbA$ is a category (because its codomain is a category). So the ``Yoneda lemma'' does not hold in general for
presheaves on $\bbA$. Therefore, those presheaves for which the Yoneda lemma happens to hold, deserve some special attention.
\begin{definition}\label{412} A {\em contravariant Yoneda presheaf} on a $\Q$-semicategory $\bbA$ is
a
$\phi\in\P{\bbA}$ satisfying, for all $a\in\bbA_0$,
\begin{equation}\label{413}
\phi(a)=\bigwedge_{a\in\bbA_0}[\bbA(a,x),\phi(x)].
\end{equation}
\end{definition} 
The hom-arrows in $\P{\bbA}$ being what they are, this definition thus says that $\phi(a)=\P{\bbA}(Y_{\bbA}a,\phi)$ for a Yoneda presheaf. In other terms, considering the $\Q$-semicategory $\bbA$ as endo-distributor on the free $\Q$-category $\ol{\bbA}$, a contravariant presheaf $\phi\:*_A\dist\ol{\bbA}$ on $\bbA$ is Yoneda if and only if
$[\bbA,\phi]=\phi$ in the quantaloid $\Dist(\Q)$. 
\par
We will denote $\Y\bbA$ for the full subcategory of $\P{\bbA}$ determined
by the Yoneda presheaves on $\bbA$.
\begin{proposition}\label{414} A $\Q$-semicategory $\bbA$ is a category if and only if all representable contravariant presheaves on $\bbA$ are Yoneda presheaves.
\end{proposition}
\proof For the non-trivial implication note that a representable
$Y_{\bbA}a$ being Yoneda implies that $\bbA(a,a)=(Y_{\bbA}a)(a)=\P{\bbA}(Y_{\bbA}a,Y_{\bbA}a)\geq 1_{ta}$, because $\P{\bbA}$ is a $\Q$-category. So when all representables are Yoneda, then $\bbA$ is in fact a category.
\endofproof 
\begin{example}\label{ex4}
Take again a strict order $(A,<)$ viewed as $\2$-semicategory $\bbA$. A presheaf 
on $\bbA$, i.e.~a downset $D\subseteq A$ of the free order $(A,\leq)$, is Yoneda if and only if $a\in D\iff \ddown a\subseteq D$.
\end{example}

\subsection*{Regular presheaves}

Every contravariant presheaf on a $\Q$-category $\bbC$ is canonicially a ``colimit of
representables'': for $\phi\in\P\bbC$, $\colim(\phi,Y_{\bbC})=\phi$. In words: the $\phi$-weighted colimit of the Yoneda embedding $Y_{\bbC}\:\bbC\to\P\bbC$ equals $\phi$.
\par
Consider now a $\Q$-semicategory $\bbA$, its
Yoneda semifunctor $Y_{\bbA}\:\bbA\to\P{\bbA}$, and a contravariant presheaf $\phi\in\P{\bbA}$. We may view $\phi$ as distributor into the free $\Q$-category $\ol{\bbA}$ and consider the $\phi$-weighted colimit of $\ol{Y_{\bbA}}\:\ol{\bbA}\to\P\bbA$:
$$\xymatrix@=15mm{\ol{\bbA}\ar[r]^{\ol{Y_{\bbA}}} & \P{\bbA} \\
\ast_X\ar[u]|{\distsign}^{\phi}\ar@{.>}[ur]_{\colim(\phi,\ol{Y_{\bbA}})}}$$
This colimit surely exists, since $\P{\bbA}$ is a cocomplete $\Q$-category ($\P\bbA=\P\ol{\bbA}$ is the free cocompletion of the free $\Q$-category $\ol{\bbA}$), but it is in
general different from (the constant functor pointing at) $\phi$! Namely, with the techniques for calculating colimits in presheaf categories it is easily seen that $\colim(\phi,\ol{Y_{\bbA}})=\bbA\tensor\phi$. Therefore the following definition is called for.
\begin{definition}\label{418} A {\em regular contravariant presheaf} on a $\Q$-semicategory $\bbA$ is a contravariant presheaf $\phi\in\P{\bbA}$ satisfying, for all $a\in\bbA_0$,
\begin{equation}\label{417}
\phi(a)=\bigvee_{x\in\bbA_0}\bbA(a,x)\circ\phi(x).
\end{equation}
\end{definition} 
That is to say, $\phi$ is a regular presheaf on $\bbA$ if and only if
$\phi=\bbA\tensor\phi$ in $\Dist(\Q)$, with $\bbA$ viewed as endo-distributor on $\ol{\bbA}$, and $\phi$ as distributor into $\ol{\bbA}$. 
\par
By $\R\bbA$ we denote the full subcategory of $\P{\bbA}$ determined
by the regular presheaves. The following must be compared with \ref{414}.
\begin{definition}\label{421} A $\Q$-semicategory $\bbA$ is {\em regular} when all representable presheaves on $\bbA$ are regular presheaves.
\end{definition}
Using (\ref{417}) this may be put as follows: $\bbA$ is a regular semicategory if and only if
\begin{equation}\label{417.1}
\bbA(a',a)=\bigvee_{x\in\bbA_0}\bbA(a',x)\circ\bbA(x,a).
\end{equation}
\par
Every $\Q$-category is a regular semicategory, but the converse is not
true. And neither is every $\Q$-semicategory regular. So the notion `regular $\Q$-semicategory' is strictly weaker than `$\Q$-category', but still {\em exactly} strong enough to allow for an interesting notion of ``presheaf which is canonically the
colimit of representables''.
\begin{example}\label{ex5}
For a presheaf on a strict order $(A,<)$, i.e.~a downset $D$ of the free order $(A,\leq)$, regularity means that, if $d\in D$ then there exists some $d'\in D$ such that $d<d'$. (In other terms, we ask that $D=\bigcup\{\ddown d\mid d\in D\}$.) So, given $a\in A$, the representable $\ddown a$ is regular if and only if for every $x<a$ there exist a $y$ such that $x<y<a$. Regularity of $(A,<)$ thus means that between every $x<z$ there lies a $y$.
\end{example}
\begin{example}\label{ex5.1}
More generally, every continous order $(A,\leq)$ with associated way-below relation $\ll$ determines a regular $\2$-category $\bbA$ by $\bbA_0=A$ and $\bbA(a,b)=1\iff a\ll b$. Indeed, way-below satisfies the ``interpolation property'' [Gierz\etal, 1980] which says that, if $x\ll z$ then $x\ll y\ll z$ for some $y$.
\end{example}
\begin{example}\label{ex5.3}
An $\Omega$-set $(A,[\cdot=\cdot])$ is, viewed as $\Omega$-semicategory $\bbA$, regular: this is due to symmetry of $\bbA(-,-)=[\cdot = \cdot]$. An $\Omega$-subset of $(A,[\cdot=\cdot])$ is precisely the same thing as a regular presheaf on $\bbA$; so $\R\bbA$ is the locale of subobjects of $\bbA$. In this framework, a representable presheaf has also been called a ``singleton''; and regularity has been interpreted as ``every subset is the union of its points''. Note that non-regular presheaves are definitely not called for in the theory of $\Omega$-sets. Take for example the $\Omega$-set $(\{*\}, [*=*]=u)$ with $u\in\Omega$ which is different from the top element. Then the (not necessarily regular) presheaves on this $\Omega$-set are formally equal to the $\Omega$-subsets of the $\Omega$-set $(\{*\}, [*=*]=\top)$. So amongst the ``parts'' of an element $*$ defined on $u\in\Omega$ would be a global element! The regularity condition on the presheaves excludes such anomalies.
\end{example}

\subsection*{Unity and identity of adjointly opposites}

In a sense, the requirement for a presheaf $\phi$ on a $\Q$-semicategory $\bbA$ to be regular is ``adjoint'' to the requirement for it to be Yoneda---compare (\ref{413}) and (\ref{417}). This can be made precise. 
\begin{theorem}\label{422} For a regular $\Q$-semicategory $\bbA$, the functors\footnote{We view a given $\Q$-semicategory $\bbA$ as endo-distributor on the free $\Q$-category $\ol{\bbA}$, and the presheaves on $\bbA$ as presheaves on $\ol{\bbA}$, so that compositions like $\bbA\tensor\psi$ and liftings like $[\bbA,\theta]$ (make sense and) are calculated in $\Dist(\Q)$.}
$$i\:\R\bbA\to\P{\bbA}\:\phi\mapsto\phi,\hspace{3mm}
j\:\P{\bbA}\to\R\bbA\:\psi\mapsto\bbA\tensor\psi,\hspace{3mm}
k\:\R\bbA\to\P{\bbA}\:\theta\mapsto[\bbA,\theta]$$
constitute adjunctions $i\dashv j\dashv k$ in $\Cat(\Q)$. Moreover, $i$ and $k$ are
fully faithful, and the image of $k$ is precisely $\Y\bbA$.
\end{theorem}
\proof  It is quite straightforward that $i$, $j$ and $k$ are well-defined functors: this depends on obvious calculations with compositions and liftings in the quantaloid $\Dist(\Q)$. As for the adjunctions:
\begin{compactitemize}
\item $i\dashv j$ means that, for all $\phi\in\R\bbA$ and all $\psi\in\P{\bbA}$, $\P{\bbA}(i\phi,\psi)=\R\bbA(\phi,j\psi)$ in $\Q$, or equivalently
$[\phi,\psi]=[\phi,\bbA\tensor\psi]$ in $\Dist(\Q)$. Using regularity of $\phi$ we have $[\phi,\bbA\tensor\psi]=[\bbA\tensor\phi,\bbA\tensor\psi]$, and reckoning that $\bbA\leq\ol{\bbA}$ it follows that $\bbA\tensor\psi\leq\psi$. So the two inequalities
$$[\phi,\psi]\leq[\bbA\tensor\phi,\bbA\tensor\psi]\mbox{ and }[\phi,\bbA\tensor\psi]\leq[\phi,\psi],$$
hold (by obvious properties of liftings in any quantaloid, thus also in $\Dist(\Q)$) and prove the claim.
\item $j\dashv k$ means that, for all $\psi\in\P{\bbA}$ and all $\theta\in\R\bbA$,
$\R\bbA(j\psi,\theta)=\P{\bbA}(\psi,k\theta)$ in $\Q$, or equivalently
$[\bbA\tensor\psi,\theta]=[\psi,[\bbA,\theta]]$ in $\Dist(\Q)$. This identity is valid in any quantaloid thus also in
$\Dist(\Q)$.
\end{compactitemize} 
The fully faithfulness of (the full inclusion) $i$ is trivial; and that
of $k$ thus follows from the adjunctions. To see that $k(\R\bbA)=\Y\bbA$, we need to show that
$\phi\in\Y\bbA$ if and only if there exists a $\psi\in\R\bbA$ such that
$\phi=k\psi$. Note first that $k\psi=[\bbA,\psi]=\R\bbA(Y_{\bbA}-,\psi)$. Now if $\phi\in\Y\bbA$ then
$$\phi= \P{\bbA}(Y_{\bbA}-,\phi)=\P\bbA(i(Y_{\bbA}-),\phi)
=\R\bbA(Y_{\bbA}-,j\phi)=k(j(\phi))$$
so $\phi=k(j(\phi))$ for $j(\phi)\in\R\bbA$.
Conversely, if $\phi=k(\psi)$ for some
$\psi\in\R\bbA$ then
$$\P{\bbA}(Y_{\bbA}-,\phi) 
=\P{\bbA}(Y_{\bbA}-,k(\psi))
=\R\bbA(j(Y_{\bbA}-),\psi)
=\R\bbA(Y_{\bbA}-,\psi)=k(\psi)
=\phi$$
so $\phi\in\Y\bbA$. We used that $j\circ Y_{\bbA}=Y_{\bbA}$, i.e.~that $\bbA$ is
regular.
\endofproof  
As a result the $\Q$-category
$\R\bbA$ is cocomplete, and colimits in $\R\bbA$ are calculated ``as in $\P\bbA$'': for a weighted diagram 
$$\xymatrix@=15mm{\bbD\ar[r]|{\distsign}^{\Theta} &
\bbC\ar[r]^F & \R\bbA}$$
the $\Theta$-weighted colimit of $F$ is
the functor
\begin{equation}\label{431.0}
\colim(\Theta,F)\:\bbD\to\R\bbA\:d\mapsto\Phi_F\tensor\Theta(-,d),
\end{equation}
where $\Phi_F\:\bbC\dist\ol{\bbA}$ has elements $\Phi_F(a,c)=(Fc)(a)$. (Indeed, this is how the colimit of $i\circ F$ would be calculated in $\P\bbA$.)  
Further, the corestriction of the fully faithful functor $k\:\R\bbA\to\P\bbA$ to its image gives an equivalence $\R\bbA\simeq\Y\bbA$ (whose inverse equivalence is a restriction of $j$). So also $\Y\bbA$ is cocomplete.

\subsection*{Alternative description of regularity} 

Regular presheaves and regular semicategories were first studied {\em as such}  in the $\V$-enriched case by [Moens\etal, 2002]. But their definition for `regular presheaf on a $\V$-semicategory' is different from -- but equivalent to -- ours (modulo a ``translation'' from $\V$-enriched structures to $\Q$-enriched structures): basically because they use the notion of ``weighted colimit of a diagram without units''.
\par
A semidistributor $\Phi\:\bbB\dist\bbA$  between $\Q$-semicategories together with a semifunctor $F\:\bbA\to\bbC$ with codomain a $\Q$-category\footnote{It is essential to consider a {\em category} $\bbC$ (whereas
$\bbA$ and $\bbB$ may be semicategories) for application of the universal property of free $\Q$-categories further on.} constitute a {\em weighted diagram without units} in the $\Q$-category $\bbC$. Due to (\ref{407}) and (\ref{408}) there is a corresponding weighted diagram ``with freely added units'', with weight $\Phi\:\ol{\bbA}\dist\ol{\bbB}$ and functor $\ol{F}\:\ol{\bbA}\to\bbC$; this diagram may or may not have a weighted colimit in the usual sense. Now, by definition, the {\em $\Phi$-weighted colimit of $F$} is the semifunctor
$\colim(\Phi,F)\:\bbB\to\bbC$ for which, by the universal property of $\ol{\bbB}$,
\begin{equation}\label{431}
\ol{\colim(\Phi,F)}=\colim(\Phi,\ol{F}).
\end{equation} 
The following proposition states the universal property of such a colimit on a weighted diagram without units,
without reference to the free categories; this was taken in [Moens\etal, 2002] as {\em definition}.
\begin{proposition}\label{432} 
Given a semidistributor $\Phi\:\bbA\dist\bbB$ between $\Q$-semicategories together with a
semifunctor $F\:\bbB\to\bbC$ into a $\Q$-category, a semifunctor $G\:\bbA\to\bbC$ is the $\Phi$-weigh\-ted colimit
of $F$ if and only if, for all $a\in\bbA_0$ and $c\in\bbC_0$, $$\bbC(Ga,c)=\bigwedge_{b\in\bbB_0}\Big[\Phi(b,c),\bbC(Fb,c)\Big],$$ the liftings on the right hand side being calculated in the base quantaloid $\Q$.
\end{proposition}
\proof Making (\ref{431}) explicit, $G$ is
the
$\Phi$-weighted colimit of $F$ if and only if $\ol{G}$ is the $\Phi$-weighted colimit of $\ol{F}$, that is to say,
$\bbC(\ol{G}-,-)=[\Phi,\bbC(\ol{F}-,-)]$ in $\Dist(\Q)$. 
As object-mappings $G$ and $\ol{G}$, resp.~$F$ and $\ol{F}$, are identical; together with the explicit formula relating liftings in $\Dist(\Q)$ to liftings in $\Q$, this gives the result.
\endofproof 
\par
For any $\Q$-semicategory $\bbA$ and $\phi\in\P{\bbA}$ we may now consider the weighted colimit of the diagram without units
$$\xymatrix@=15mm{
\bbA\ar[r]^{Y_{\bbA}} & \P\bbA \\
\ast_X\ar[u]^{\phi}|{\distsign}\ar@{.>}[ur]_{\colim(\phi,Y_{\bbA})}}$$
Trivially the presheaf $\phi$ is regular in the sense of \ref{418} -- i.e.~equal to $\phi$-weighted colimit of
$\ol{Y_{\bbA}}\:\ol{\bbA}\to\P\bbA$ -- if and only if it is
the $\phi$-weighted colimit of $Y_{\bbA}\:\bbA\to\P\bbA$. By \ref{432} we get the following characterization in terms of liftings in $\Q$, which was taken as {\em definition} of `regular presheaf on a semicategory' in [Moens\etal, 2002].
\begin{proposition}\label{433} A presheaf $\phi\in\P{\bbA}$ on a
$\Q$-semicategory $\bbA$ is regular if and only if, for every $\psi\in\P{\bbA}$, $$\P{\bbA}(\phi,\psi)=\bigwedge_{a\in\bbA_0}\Big[\phi(a),\P\bbA(Y_{\bbA}a,\psi)\Big].$$ 
\end{proposition}

\subsection*{Covariant presheaves}

Dualizing the theory of {\em contravariant} presheaves gives the theory of {\em covariant} presheaves. Recalling that, for a $\Q$-category $\bbC$, the $\Q$-category of covariant presheaves is denoted $\P\+\bbC$, we denote the {\em $\Q$-category of covariant presheaves on a $\Q$-semicategory $\bbA$} as $\P\+\bbA$. It is straightforward to define Yoneda and regular covariant presheaves; the respective full subcategories of $\P\+\bbA$ are naturally denoted as $\Y\+\bbA$ and $\R\+\bbA$. A $\Q$-semicategory $\bbA$ is a $\Q$-category if and only if all representable covariant presheaves are Yoneda (compare with \ref{414}); and $\bbA$ is a regular $\Q$-semicategory if and only if all representable covariant presheaves are regular (compare with \ref{421}). That is to say, it comes to the same thing to ask for either the {\em contravariant} or the {\em covariant} presheaves to be Yoneda, resp.~regular.
\begin{example}\label{ex5.2}
Let $(A,\leq)$ be a continuous order; denote $\bbA$ for the regular $\2$-semicategory whose hom-arrows classify the way-below relation. A covariant regular presheaf on $\bbA$ is the same thing as a Scott-open subset of $A$, and a contravariant Yoneda presheaf is the same thing as a Scott-closed subset; this requires some straightforward calculations with the way-below relation. In particular, the $\2$-category $\R\+\bbA$ is the Scott-topology on $\bbA$ (i.e.~the collection of Scott-opens of $(A,\leq)$ ordered by inclusion).
\end{example}

\section{Regular semidistributors, Morita equivalence}\label{465}

\subsection*{Calculus of regular semidistributors}

Generalizing (\ref{417}) and (\ref{417.1}) we define a `regular semidistributor' as follows.
\begin{definition}\label{466}
A semidistributor $\Phi\:\bbA\dist\bbB$ between $\Q$-semicategories is {\em regular} if, for all $a\in\bbA_0$ and $b\in\bbB_0$,
$$\bigvee_{x\in\bbA_0}\Phi(b,x)\circ\bbA(x,a)=\Phi(b,a)=\bigvee_{y\in\bbB_0}\bbB(b,y)\circ\Phi(y,a).$$
\end{definition}
When viewing the semidistributor $\Phi\:\bbA\dist\bbB$ as distributor $\Phi\:\ol{\bbA}\dist\ol{\bbB}$, and both $\bbA$ and $\bbB$ as endo-distributors on the respective free categories $\ol{\bbA}$ and $\ol{\bbB}$, this definition says that $\bbB\tensor\Phi=\Phi=\Phi\tensor\bbA$ in $\Dist(\Q)$.
\par
It follows straightforwardly that, for a $\Q$-semicategory $\bbA$, 
\begin{compactitemize}
\item a semidistributor $\phi\:*_X\dist\bbA$ is a regular presheaf on $\bbA$ (in the sense of
\ref{418}) if and only if it is a regular semidistributor (in the sense
of \ref{466}),
\item $\bbA$ is a regular $\Q$-semicategory (in the sense of \ref{421}) if and only if
the semidistributor $\bbA\:\bbA\dist\bbA$ is regular (in the sense of \ref{466}).
\end{compactitemize}
So there can be no confusion when we use the word `regular'!
\par
As mentioned before, $\Q$-semicategories and semidistributors do not form a quantaloid, the problem being that ``a semicategory $\bbA$ is {\em not} the identity semidistributor on itself''. The notions of {\em regular} semicategory and {\em regular} semidistributor solve that problem.
\begin{proposition}\label{466.2}
Regular $\Q$-semicategories are the objects, and regular semidistributors the arrows,
of a quantaloid
$\RSDist(\Q)$ in which local suprema are as in (\ref{408.1}), composition is as in (\ref{408.0}), and identities are as in (\ref{408.2}). It contains $\Dist(\Q)$ as full subquantaloid.
\end{proposition}
\proof
It is straightforward that $\RSDist(\Q)$ is a quantaloid: the condition in \ref{466} says that, for a {\em regular} semidistibutor $\Phi\:\bbA\dist\bbB$ between {\em regular} $\Q$-semicategories, the endo-semidistributors $\bbA\:\bbA\dist\bbA$ and $\bbB\:\bbB\dist\bbB$ are {\em identities} for the composition law in (\ref{408.0}). That $\Dist(\Q)$ is a full subquantaloid, follows from previous remarks, plus the fact that composition of distributors between categories is not different from composition of (regular) semidistributors between (regular) semicategories.
\endofproof
The appendix contains a more elegant argument: $\RSDist(\Q)$ is the universal split-idempotent completion of $\Matr(\Q)$ in $\QUANT$, that is to say, it is the Cauchy completion of $\Q$ in $\QUANT$.
\begin{example}\label{ex6.0}
Viewing $\Omega$-sets $(A,[\cdot=\cdot])$ and $(B,[\cdot=\cdot])$ as regular $\Omega$-semicategories $\bbA$ and $\bbB$, a morphism of $\Omega$-sets is precisely a left adjoint regular semidistributor. That is to say, the topos of sheaves on a locale $\Omega$ may be described (up to equivalence) as: objects are symmetric regular $\Omega$-semicategories, and arrows are left adjoint regular semidistributors.
\end{example}
\par
Given the formula for liftings in $\RSDist(\Q)$ (see the appendix), we can give a useful alternative description of (the hom-arrows in) the category of regular presheaves on a {\em regular} $\Q$-semicategory.
\begin{corollary}\label{466.3}
For a regular $\Q$-semicategory $\bbA$, the category $\R\bbA$ of regular contravariant presheaves may be described as:
\begin{compactitemize}
\item objects: an object of type $X\in\Q$ is a regular semidistributor
$\phi\:*_X\dist\bbA$,
\item hom-arrows: for regular semidistributors $\phi\:*_X\dist\bbA$ and $\psi\:*_Y\dist\bbA$, $\R\bbA(\psi,\phi)$ is (the single element of) the lifting $[\psi,\phi]$ in the quantaloid $\RSDist(\Q)$.
\end{compactitemize}
\end{corollary}
\proof
Let us, for the sake of clarity, write $\R_1\bbA$ for the $\Q$-category defined in the statement of this proposition, and $\R_2\bbA$ for the $\Q$-category of regular presheaves on $\bbA$, as defined previously (as full subcategory of $\P\bbA$). Then, by the remarks following \ref{466}, $(\R_1\bbA)_0=(\R_2\bbA)_0$ as $\Q$-typed object-sets. That the hom-arrows coincide too, follows from (\ref{app4}) in the appendix: let $\phi$ and $\psi$ denote two regular presheaves on $\bbA$, so we may either view them as distributors $\phi\:*_X\dist\ol{\bbA}$ and $\psi\:*_Y\dist\ol{\bbB}$ satisfying $\bbA\tensor\phi=\phi$ and $\bbA\tensor\psi=\psi$, or we may view them as regular semidistributors $\phi\:*_X\dist\bbA$ and $\psi\:*_Y\dist\bbA$. Since $*_X$ and $*_Y$ are {\em categories} we have that
$$\R_1\bbA(\psi,\phi)=[\psi,\phi]^{\RSDist(\Q)}=[\psi,\phi]^{\Dist(\Q)}=\R_2\bbA(\psi,\phi).$$
So $\R_1\bbA=\R_2\bbA$ as claimed.
\endofproof

\subsection*{An aspect of Morita equivalence}

Consider a regular semidistributor between regular $\Q$-semicategories, like $\Phi\:\bbA\dist\bbB$. To the associated distributor $\Phi\:\ol{\bbA}\dist\ol{\bbB}$ corresponds a functor $F_{\Phi}\:\ol{\bbA}\to\P\bbB$ (by the classifying property of $\P\bbB=\P\ol{\bbB}$) which actually lands in $\R\bbB$: every $F_{\Phi}(a)=\Phi(-,a)$ is a regular presheaf on $\bbB$. Since $\R\bbB$ is a cocomplete $\Q$-category, the left Kan extension of $F_{\Phi}\:\ol{\bbA}\to\R\bbB$ along $\ol{Y_{\bbA}}\:\ol{\bbA}\to\R\bbA$, denoted $\<F_{\Phi},\ol{Y_{\bbA}}\>\:\R\bbA\to\R\bbB$, exists and is pointwise; so $\<F_{\Phi},\ol{Y_{\bbA}}\>(\phi)=\colim(\phi,F_{\Phi})$. But with (\ref{431.0}), and using the description of $\R\bbA$ and $\R\bbB$ as in \ref{466.3}, a practical explicit description of $\<F_{\Phi},\ol{Y_{\bbA}}\>$ can be given:
$$\<F_{\Phi},\ol{Y_{\bbA}}\>\:\R\bbA\to\R\bbB\:\phi\mapsto\Phi\tensor\phi,$$
the right hand side being simply a composition of arrows in $\RSDist(\Q)$.
The following proposition is then really a generalization of what is known for distributors between $\Q$-categories.
\begin{proposition}\label{470}
For any quantaloid $\Q$, 
\begin{equation}\label{471}
\RSDist(\Q)\to\Cocont(\Q)\:\Big(\Phi\:\bbA\dist\bbB\Big)\mapsto
\Big(\Phi\tensor-\:\R\bbA\to\R\bbB\Big)
\end{equation}
is a 2-functor which is locally an equivalence.
\end{proposition}
\proof
First observe that, for any regular semidistributor $\Phi\:\bbA\dist\bbB$ between regular semicategories, the functor $\R\bbA\to\R\bbB\:\theta\mapsto\Phi\tensor\theta$ has a right adjoint in $\Cat(\Q)$, namely $\R\bbB\to\R\bbA\:\psi\mapsto[\Phi,\psi]$ (the lifting being calculated in $\RSDist(\Q)$). Since furthermore the action in (\ref{471}) is determined by composition of arrows in $\RSDist(\Q)$, it is then obviously a well-defined 2-functor. This 2-functor (between locally ordered categories) reflects the local order: if $\Phi,\Psi\:\bbA\bidist\bbB$ in $\RSDist(\Q)$ satisfy 
$\Phi\tensor\theta\leq\Psi\tensor\theta$ for all $\theta\in\R\bbA$
then in particular 
$$\Phi(-,a)=\Phi\tensor\bbA(-,a)\leq\Psi\tensor\bbA(-,a)=\Psi(-,a)$$ 
for any $a\in\bbA_0$ so $\Phi\leq\Psi$ in $\RSDist(\Q)$. Finally we must see that (\ref{471}) is locally surjective. Given a cocontinuous functor $F\:\R\bbA\to\R\bbB$, putting 
$$\Phi(b,a)=F(Y_{\bbA}a)(b)\mbox{ for $a\in\bbA_0$, $b\in\bbB_0$}$$
defines a regular semidistributor from $\bbA$ to $\bbB$: on the one hand $\bbB\tensor\Phi=\Phi$ because every $\Phi(-,a)=F(Y_{\bbA}a)$ is a regular presheaf on $\bbB$, so $\Phi$ is regular in $\bbB$; on the other hand 
$$\Phi\tensor\bbA(-,a)=\colim(Y_{\bbA}a,F\circ Y_{\bbA})=F\circ\colim(Y_{\bbA}a,Y_{\bbA})=F(Y_{\bbA}a)=\Phi(-,a),$$
so $\Phi$ is also regular in $\bbA$. (Note the use of regularity of $\bbA$, and the calculation of colimits in $\R\bbB$ cf.~(\ref{431.0})!)
This regular semidistributor $\Phi\:\bbA\dist\bbB$ satisfies, for $\theta\in\R\bbA$, 
$$\Phi\tensor\theta=\colim(\theta,F\circ Y_{\bbA})=F\circ\colim(\theta,Y_{\bbA})=F(\theta)$$
so $F$ is the image of $\Phi$ under (\ref{471}).
\endofproof
Recalling that cocomplete $\Q$-categories are equivalent in $\Cat(\Q)$ if and only if they are equivalent in $\Cocont(\Q)$, and recalling that $\R\bbA$ and $\Y\bbA$ are equivalent, we may record the following corollary (where the dual to \ref{422} and \ref{470} are assumed).
\begin{corollary}\label{472}
For regular $\Q$-semicategories $\bbA$ and $\bbB$, the following are equivalent:
\begin{enumerate}
\item $\bbA\cong\bbB$ in $\RSDist(\Q)$,
\item $\R\bbA\simeq\R\bbB$ in $\Cat(\Q)$,
\item $\R\+\bbA\simeq\R\+\bbB$ in $\Cat(\Q)$,
\item $\Y\bbA\simeq\Y\bbB$ in $\Cat(\Q)$,
\item $\Y\+\bbA\simeq\Y\+\bbB$ in $\Cat(\Q)$.
\end{enumerate}
\end{corollary}
When in the above $\bbA$ and $\bbB$ happen to be $\Q$-categories, then $\R\bbA=\Y\bbA=\P\bbA$ and $\R\+\bbA=\Y\+\bbA=\P\+\bbA$, and likewise for $\bbB$, so that the stated equivalence is an aspect of the Morita equivalence of $\Q$-categories. Therefore,
two regular $\Q$-semicategories $\bbA$ and $\bbB$ can be called {\em Morita equivalent} if the equivalent conditions in \ref{472} hold\footnote{Morita equivalence for $\Q$-{\em categories} also means that the categories have equivalent {\em Cauchy completions}; but this result is simply impossible in the theory of regular $\Q$-semicategories: in general such semicategories {\em do not have a Cauchy completion} (in a suitable sense of the word). For a study of those $\Q$-semicategories that do admit a Cauchy completion, we refer to [Stubbe, 2004b], for the matter is too involved to include it here.}.
In general it is not true that every regular $\Q$-semicategory is Morita equivalent to a $\Q$-category, which makes the theory of regular $\Q$-semicategories a subject that is really different from the theory of $\Q$-categories.
\begin{example}\label{ex6}
Let $\bf 3$ denote the three-element chain $0<e<1$ viewed as quantaloid with one object. Let $\bbA$ be the $\bf 3$-semicategory with one object $*$ and hom-arrow $\bbA(*,*)=e$: since $e\wedge e=e$, this semicategory is regular. The $\bf 3$-category $\R\bbA$ of regular contravariant presheaves has two (non-isomorphic) objects, $0$ and $e$. The (only) $\bf 3$-category $\bbC$ with one object necessarily has $\bbC(*,*)=1$, and there are three (non-isomorphic) presheaves on $\bbC$: $0$, $e$ and $1$. There can be no equivalence between $\R\bbA$ and $\P\bbC$, so $\bbA$ and $\bbC$ are not Morita equivalent. Considering a $\bf 3$-category with more objects only increases the number of presheaves on it, and makes it even more impossible for such a category to be Morita equivalent to $\bbA$.
\end{example}

\section{Regular semifunctors}

Any functor between $\Q$-categories induces an
adjoint pair of distributors; this is a simple but crucial fact for the development of $\Q$-category theory. For $\Q$-semi\-cat\-e\-go\-ries things are a little bit more involved: only ``regular semifunctors'' between regular semicategories do the trick.

\subsection*{Inducing adjoints}
 
A semifunctor $F\:\bbA\to\bbB$ between regular $\Q$-semicategories determines semidistributors $\bbB(-,F-)\:\bbA\dist\bbB$ and $\bbB(F-,-)\:\bbB\dist\bbA$, whose distributor-arrows are respectively
$$\bbB(b,Fa)\:ta\to tb\mbox{ and }\bbB(Fa,b)\:tb\to ta\mbox{, for $a\in\bbA_0$ and $b\in\bbB_0$.}$$
But these semidistributors needn't be regular! And when they are not regular, they are not arrows in a suitable bicategory -- {\it in casu} $\RSDist(\Q)$ -- and so they cannot be ``adjoint''! The following definition is thus called for.
\begin{definition}\label{485}
A semifunctor $F\:\bbA\to\bbB$ between $\Q$-semicategories is {\em regular} when both semidistributors 
$\bbB(-,F-)\:\bbA\dist\bbB$ and $\bbB(F-,-)\:\bbB\dist\bbA$ are regular.
\end{definition}
If $\bbA$ and $\bbB$ are {\em regular} $\Q$-semicategories, then $\bbB(-,F-)$ and $\bbB(F-,-)$ are arrows in $\RSDist(\Q)$.
\begin{proposition}\label{487} Regular $\Q$-semicategories are the objects, 
and regular semifunctors the arrows, of a (non-full) subcategory $\RSCat(\Q)$ of $\SCat(\Q)$. It contains $\Cat(\Q)$ as full subcategory.
\end{proposition} 
\proof
The identity semifunctor $1_{\bbA}\:\bbA\to\bbA$ on a $\Q$-semicategory is regular
if and only if the semicategory itself is regular. If $F\:\bbA\to\bbB$ and $G\:\bbB\to\bbC$ are regular semifunctors between regular
$\Q$-semicategories, then $\bbB(-,F-)\:\bbA\dist\bbB$ and $\bbC(-,G-)\:\bbB\dist\bbC$ are
arrows in $\RSDist(\Q)$. Composition of these $\RSDist(\Q)$-arrows gives $\bbC(-,G\circ F-)$, with composition of semifunctors in $\SCat(\Q)$---so the latter is a regular semidistributor. Similar for $\bbC(G\circ F-, -)$. So $G\circ F$, computed in $\SCat(\Q)$, is a regular semifunctor. This proves that $\RSCat(\Q)$ is a subcategory of $\SCat(\Q)$; that $\Cat(\Q)$ is a full subcategory of $\RSCat(\Q)$ is straightforward.
\endofproof 
\begin{proposition}\label{488} A regular semifunctor 
$F\:\bbA\to\bbB$ between regular $\Q$-semi\-cat\-e\-go\-ries determines an adjoint pair $\bbB(-,F-)\dashv\bbB(F-,-)$ in $\RSDist(\Q)$. Actually, 
\begin{equation}\label{488.1}
\RSCat(\Q)\to\RSDist(\Q)\:\Big(F\:\bbA\to\bbB\Big)\mapsto\Big(\bbB(-,F-)\:\bbA
\dist\bbB\Big)
\end{equation}
is functorial. (And $F\mapsto\bbB(F-,-)$ is contravariantly functorial.)
\end{proposition}
\proof  
The adjunction's unit inequality $\bbB(F-,-)\tensor\bbB(-,F-)=\bbB(F-,F-)\geq\bbA$ holds by regularity of $\bbB$ and the action-inequality for $F$; the co-unit inequality 
$\bbB(-,F-)\tensor\bbB(F-,-)\leq\bbB(-,-)\tensor\bbB(-,-)=\bbB$  follows from $F(\bbA_0)\subseteq\bbB_0$ -- using the explicit expression for composition in $\RSDist(\Q)$, cf.~(\ref{408.0}) -- and once again regularity of $\bbB$.
\par
For the functoriality of (\ref{488.1}), we must verify that $1_{\bbA}\:\bbA\to\bbA$ is sent to
$\bbA\:\bbA\dist\bbA$, and for $F\:\bbA\to\bbB$ and $G\:\bbB\to\bbC$, that 
$\bbC(-,G\circ F-)=\bbC(-,G-)\tensor\bbB(-,F-)$. But this was already contained in the proof of \ref{487}, so we need not repeat it here.
\endofproof  
$\RSCat(\Q)$ now inherits local structure
from
$\RSDist(\Q)$. This makes $\RSCat(\Q)$ a locally ordered 2-category, and the functor in (\ref{488.1}) is then a 2-functor. Moreover the diagram
$$\xymatrix@=15mm{
\RSCat(\Q)\ar[r] & \RSDist(\Q) \\
\Cat(\Q)\ar[r]\ar[u] & \Dist(\Q)\ar[u]}$$
of 2-categories and 2-functors commutes: the horizontal arrows say that ``every regular semifunctor (functor) induces an adjoint pair of regular semidistributors (distributors)'', and the vertical ones are full embeddings.

\subsection*{Inadequacy of regular semifunctors}

For the theory of regular $\Q$-semicategories, \ref{488} is of much importance: it will be the starting point when we look into a theory of ``Cauchy completable'' semicategories [Stubbe, 2004b]. But some (counter)examples indicate that not all interesting ``morphisms between (regular) $\Q$-semicategories'' are regular semifunctors.
\begin{example}\label{ex7}
For the regular $\bf 3$-semicategory $\bbA$ described in \ref{ex6}, the Yoneda semifunctor $Y_{\bbA}\:\bbA\to\R\bbA$ is not regular: neither $\R\bbA(-,Y_{\bbA}-)$ nor $\R\bbA(Y_{\bbA}-,-)$ is a regular semidistributor. 
And neither is the semifunctor $i_{\bbA}\:\bbA\to\ol{\bbA}$ into its free category is regular; again neither $\ol{\bbA}(-,i_{\bbA}-)$ nor $\ol{\bbA}(i_{\bbA}-,-)$ is regular.
\end{example}
\begin{example}\label{ex9}
Let $(A,\leq)$ and $(B,\leq)$ be continuous orders; one of the many equivalent ways of saying that a map $f\:A\to B$ is Scott-continuous, is to ask that $b\ll fa$ in $B$ if and only if there exists an $x\in A$ such that $b\ll fx$ and $x\ll a$. Taking $\bbA$ and $\bbB$ to be the regular $\2$-semicategories associated with the way-below relations on $A$ and $B$, this says that $\bbA_0\to\bbB_0\:a\mapsto fa$ is an object-mapping for which $\bbB(-,f-)$ is a regular distributor from $\bbA$ to $\bbB$. But it doesn't say that this object-mapping is a semifunctor (i.e.~$a\ll a'$ doesn't imply $fa\ll fa'$, even though $x\leq y$ does imply $fx\leq fy$), and not does it say that $\bbB(f-,-)$ is a regular distributor!
\end{example}

\section{Appendix: Calculus of regular semidistributors}

\subsection*{Splitting idempotents}

For any category $\C$, the following defines a new
category
$\Idm(\C)$ of ``idempotents in $\C$'' [Freyd, 1964]:
\begin{compactitemize}
\item objects: are all the idempotents $A\Endoar{e}$, $B\Endoar{f}$,
$C\Endoar{g}$, ... in $\C$;
\item arrows: an $\Idm(\C)$-arrow $b\:e\dist f$ between idempotents $A\Endoar{e}$ and $B\Endoar{f}$ is a $\C$-arrow $b\:A\to B$ such that $b\circ e=b=f\circ b$;
\item composition: the composition of $b\:e\dist f$ and $c\:f\dist g$ is written $c\tensor b\:e\dist g$ computed as
$c\tensor b=c\circ b$ in $\C$;
\item identity: on an idempotent $ A\Endoar{e}$, $e\colon e\dist e$ itself is the identity.
\end{compactitemize} 
The category $\C$ can be embedded in $\Idm(\C)$ by sending $f\:A\to B$ to $f\:1_A\dist 1_B$; let us write this as $k_{\C}\:\C\to\Idm(\C)$. This embedding is full, and it is an equivalence if and only if all idempotents in $\C$ split. All idempotents in $\Idm(\C)$ split (so ``taking idempotents'' is an idempotent procedure); in fact, $k_{\C}$ is the universal split-idempotent completion of $\C$ in $\CAT$. Further
it is straightforward that, if a category $\C$ admits (co)products, then so does the category 
$\Idm(\C)$.
\par
Taking $\C$ to be (the underlying category of) a (possibly large) quantaloid $\Q$ in the above, it is quite obvious that $\Idm(\Q)$ is a quantaloid too: arrows in $\Idm(\Q)$ are ordered as they are ordered in $\Q$, and so the embedding $k_{\Q}\:\Q\to\Idm(\Q)$ is a quantaloid homomorphism. Precisely, $\Q\to\Idm(\Q)$ is the universal split-idempotent completion of $\Q$ in $\QUANT$. ([Koslowski, 1997] studies the splitting of ``idempolads'' in more general bicategories, of which $\Idm(\Q)$ is a particularly simple example.)
\par
As any quantaloid, $\Idm(\Q)$ is closed. For example, for $b\:e\dist f$ and $c\:g\dist f$, the lifting of $b$ through $c$ in $\Idm(\Q)$ can be calculated as
$$[c,b]^{\Idm(\Q)}=g\circ[c,b]^{\Q}\circ e,$$
the composition and lifting on the right hand being calculated in $\Q$. (One verifies that the expression on the right hand equals $\bigvee\{d\:e\dist g\mbox{ in }\Idm(\Q)\mid c\tensor d\leq b\}$.) There is a similar formula for extensions. 

\subsection*{Adding direct sums}

Two other universal constructions on a quantaloid $\Q$ play an important r\^ole: the quantaloid $\Matr(\Q)$ of ``matrices with elements in $\Q$'', for which $i_{\Q}\:\Q\to\Matr(\Q)$ is the direct-sum completion of $\Q$ in $\QUANT$, and the quantaloid $\Mnd(\Q)$ of ``monads and bimodules in $\Q$'', for which $j_{\Q}\:\Q\to\Mnd(\Q)$ is the split-monad completion of $\Q$ in $\QUANT$. For details we refer to [Carboni\etal, 1987], or the appendix of [Stubbe, 2004a].
\par
It is then quite obvious that $\Mnd(\Q)$ is fully embedded in $\Idm(\Q)$ (``all monads in a quantaloid are idempotent''), and if $k_{\Q}\:\Q\to\Idm(\Q)$ is an equivalence then so is $j_{\Q}\:\Q\to\Mnd(\Q)$ (``if all idempotents split, then surely monads split''). Further,   $\Idm(\Matr(\Q))$ has all direct sums, because $\Matr(\Q)$ does so.
\par
It is but an observation that, precisely like $\Dist(\Q)=\Mnd(\Matr(\Q))$, we have $\RSDist(\Q)=\Idm(\Matr(\Q))$. In words: regular $\Q$-semicategories are exactly idempotent matrices with elements in $\Q$, and regular distributors between such semicategories are exactly matrices that are compatible with those idempotent matrices; this is really the content of (\ref{417.1}) and \ref{466}. This proves thus that $\RSDist(\Q)$ is a quantaloid in which $\Dist(\Q)$ is fully embedded. And it says that the full embedding $\Q\to\RSDist(\Q)$, sending an arrow $f\:A\to B$ onto $(f)\:*_A\dist *_B$, is the universal direct-sum-and-split-idempotent completion of $\Q$ in $\QUANT$.
\par
Note in particular that for regular semidistributors $\Phi\:\bbA\dist\bbB$ and $\Psi\:\bbC\dist\bbB$, the lifting of $\Phi$ through $\Psi$ may be calculated as
$$[\Psi,\Phi]^{\RSDist(\Q)}=\bbC\circ[\Psi,\Phi]^{\Matr(\Q)}\circ\bbA,$$
the right hand side being calculated in $\Matr(\Q)$. But, regarding $\Phi$ and $\Psi$ as distributors between free $\Q$-categories, like $\Phi\:\ol{\bbA}\dist\ol{\bbB}$ and $\Psi\:\ol{\bbC}\dist\ol{\bbB}$, also
$$\bbC\circ[\Psi,\Phi]^{\Matr(\Q)}\circ\bbA=\bbC\tensor[\Psi,\Phi]^{\Dist(\Q)}\tensor\bbA,$$ 
holds: because lifting and composition in $\Dist(\Q)$ are calculated ``as in $\Matr(\Q)$''! So, given regular semidistributors $\Phi\:\bbA\dist\bbB$ and $\Psi\:\bbC\dist\bbB$ we know that
$$[\Psi,\Phi]^{\RSDist(\Q)}=\bbC\tensor[\Psi,\Phi]^{\Dist(\Q)}\tensor\bbA,$$
where on the right hand $\Phi$ and $\Psi$ are regarded as distributors between categories, $\Phi\:\ol{\bbA}\dist\ol{\bbB}$ and $\Psi\:\ol{\bbC}\dist\ol{\bbB}$. If now $\bbA$ and $\bbC$ happen to be $\Q$-categories, then $\bbA=\ol{\bbA}$ and $\bbC=\ol{\bbC}$ so the above reduces to
\begin{equation}\label{app4}
[\Psi,\Phi]^{\RSDist(\Q)}=[\Psi,\Phi]^{\Dist(\Q)}.
\end{equation}
This is useful in the proof of \ref{466.3}.

\subsection*{Cauchy completion of a quantaloid}

It is well known that, for a small category $\C$, $k_{\C}\:\C\to\Idm(\C)$ is the Cauchy completion of $\C$ in $\CAT$, meaning (amongst other things) that it is the universal completion for absolute (co)limits, i.e.~those (co)limits that are preserved by any functor. 
\par
For a small quantaloid $\Q$, the split-idempotent completion
$k_{\Q}\:\Q\to\Idm(\Q)$ is {\em not} its
Cauchy completion in $\QUANT$, since $\Idm(\Q)$ is not Cauchy complete in $\QUANT$: surely
the splitting of idempotents is a necessary condition (any such splitting is a
(co)equalizer preserved by any quantaloid homomorphism), but -- for example -- also direct sums in a quantaloid are absolute, so a Cauchy complete quantaloid must necessarily admit these too.
But note that the quantaloid $\Idm(\Matr(\Q))$ has all direct sums and that all of its
idempotents split, so that the embedding
$\Q\to\Idm(\Matr(\Q))$ is a universal construction in $\QUANT$ adding to $\Q$ all absolute (co)equalizers and all absolute (co)products. Indeed, it is the Cauchy completion of $\Q$ in
$\QUANT$.
Although this presentation is slightly 
different from that in [Van der Plancke, 1997], a proof is given by
the elementary calculations on pp.~67--71 of that work.
\par
All this implies that, for a small quantaloid $\Q$, $\Q\to\RSDist(\Q)$ is its Cauchy completion in $\QUANT$.

\par

\end{document}